%% file: CRisp_enki.tex
\documentclass[10pt,reqno]{amsart}
\usepackage{amssymb,amscd,amsbsy}
\input{CRinput}

\begin{document}

\newcommand{\vse}{\vspace{.2in}}

\title{Functions of perturbed tuples of self-adjoint operators}


\maketitle

\vspace*{-.5cm}

\begin{center}
\Large
Fedor Nazarov$^{\rm a}$, Vladimir Peller$^{\rm b}$
\end{center}

\begin{center}
\footnotesize
{\it$^{\rm a}$Department of Mathematics, Kent State University, Kent, OH 44242, USA\\
$^{\rm b}$Department of Mathematics, Michigan State University, East Lansing, MI 48824, USA}
\end{center}

\newcommand\fM{\frak M}
\newcommand\dg{\frak D}

\footnotesize

\noindent
{\bf Abstract.} We generalize earlier results of \cite{AP2}, \cite{AP3}, \cite{APPS2}, 
\cite{Pe1}, \cite{Pe2} to the case of functions of $n$-tuples of commuting self-adjoint 
operators. In particular, we prove that if a function $f$ belongs to the Besov space $B_{\be,1}^1(\R^n)$, then $f$ 
is operator Lipschitz and we show that if $f$ satisfies a H\"older condition of order $\a$,
then $\|f(A_1\cdots,A_n)-f(B_1,\cdots,B_n)\|\le\const\max_{1\le j\le n}\|A_j-B_j\|^\a$ for all $n$-tuples of commuting self-adjoint operators $(A_1,\cdots,A_n)$ and $(B_1,\cdots,B_n)$. We also consider the case of arbitrary moduli of continuity and the case when the operators $A_j-B_j$ belong to the Schatten--von Neumann class $\bS_p$.

\medskip

\begin{center}
{\bf\large Fonctions d'uplets d'op\'erateurs autoadjoints perturb\'es.} 
\end{center}

\medskip

\noindent
{\bf R\'esum\'e.} Dans cette note nous g\'en\'eralisons des r\'esultats de \cite{AP2}, \cite{AP3}, \cite{APPS2}, \cite{Pe1}, \cite{Pe2} en cas de fonctions d'op\'erateurs auto-adjoints et d'op\'erateurs normaux. Nous consid\'erons le probl\`eme similaire pour les fonctions de $n$-uplets d'op\'erateurs auto-adjoints qui commutent. En particulier, nous d\'emontrons que si $f$ est une fonction de la classe de Besov $B_{\be,1}^1(\R^n)$, alors
elle est lipschitzienne op\'eratorielle. En outre, nous montrons que si $f$ appartient \`a
l'espace de H\"older d'ordre $\a$, alors $\|f(A_1\cdots,A_n)-f(B_1,\cdots,B_n)\|\le\const\max_{1\le j\le n}\|A_j-B_j\|^\a$ por tous $n$-uplets $(A_1,\cdots,A_n)$ et 
$(B_1,\cdots,B_n)$ d'op\'erateurs auto-adjoints qui commutent. Nous consid\'erons aussi le cas de module de continuit\'e arbitraire et le cas o\`u les op\'erateurs $A_j-B_j$ appartiennent \`a l'espace de Schatten--von Neumann $\bS_p$.

\normalsize

\medskip

\begin{center}
{\bf\large Version fran\c caise abr\'eg\'ee}
\end{center}

\medskip

Il est bien connu (voir \cite{F}) qu'il y a des fonctions $f$ lipschitziennes sur la droite r\'eelle $\R$ qui ne sont pas {\it lipschitziennes op\'eratorielles},
c'est-\`a-dire la condition
$|f(x)-f(y)|\le\const|x-y|$,  $x,\,y\in\R,$
n'implique pas que pour tous les op\'erateurs auto-adjoints $A$ et $B$ l'in\'egalit\'e
$$
\|f(A)-f(B)\|\le\const\|A-B\|.
$$
soit vraie.
Dans \cite{Pe1} et \cite{Pe2} des conditions n\'ecessaires et des conditions suffisantes sont donn\'ees pour qu'une fonction $f$ soit lipschitzienne
op\'eratorielle. En particulier, il est d\'emontr\'e dans \cite{Pe1}
que pour qu'une fonction $f$ soit lipschitzienne op\'eratorielle il est n\'ecessaire que $f$ appartienne localement \`a l'espace de Besov $B_{11}^1(\R)$.
Cela implique aussi qu'une fonction lipschitzienne n'est pas n\'ecessairement
lipschitzienne op\'eratorielle. D'autre part, il est d\'emontr\'e dans \cite{Pe1} et \cite{Pe2} que si $f$ appartient \`a l'espace de Besov $B_{\be1}^1(\R)$, alors la fonction $f$ est lipschitzienne
op\'eratorielle.

Il se trouve que la situation change dramatiquement si l'on consid\`ere les fonctions de la classe $\L_\a(\R)$ de H\"older d'ordre $\a$, $0<\a<1$.
Il est d\'emontr\'e dans \cite{AP1} et \cite{AP2}  que si
$f\in\L_\a(\R)$, $0<\a<1$,
(c'est-\`a-dire $|f(x)-f(y)|\le\const|x-y|^\a$), alors $f$ doit \^etre
{\it h\"olderienne op\'eratorielle d'ordre} $\a$, c'est-\`a-dire
$$
\|f(A)-f(B)\|\le\const\|A-B\|^\a.
$$
Ce r\'esultat \'etait g\'en\'eralis\'e dans \cite{AP2} pour les modules de continuit\'e arbitraires. 

Les r\'esultats ci-dessus ont \'et\'e g\'eneralis\'es dans \cite{APPS1} et \cite{APPS2} au cas de fonctions d'op\'erateurs normaux. 

Dans cette note nous consid\'erons le cas de fonctions de $n$-uplets d'op\'erateurs auto-adjoints qui commutent. Il se trouve que les m\'ethodes du travail \cite{APPS2} ne marchent pas dans cette situation. Supposons que $f$ est une fonction born\'ee sur $\R^3$ dont la transform\'ee de Fourier a un support compact. On peut montrer que comme dans le cas d'op\'erateurs normaux, les fonction $\dg_1f$ et $\dg_3f$ sur $\R^3\times\R^3$ d\'efinies par
$$
(\dg_1f)(x,y)=\frac{f(x_1,x_2,x_3)-f(y_1,x_2,x_3)}{x_1-y_1},\quad
(\dg_3f)(x,y)=\frac{f(y_1,y_2,x_3)-f(y_1,y_2,y_3)}{x_3-y_3}
$$
sont des multiplicateurs de Schur (voir \S\,3 pour la d\'efinition). Toutefois, contrairement au cas $n=2$, la fonction $\dg_2$ d\'efinie par
$$
(\dg_2f)(x,y)=\frac{f(y_1,x_2,x_3)-f(y_1,y_2,x_3)}{x_2-y_2},
$$
n'est pas un multiplicateur de Schur.

Cependant, nous d\'emontrons le r\'esultat suivant:

{\it Soient $\s>0$ et $f$ une fonction dans $L^\be(\R^n)$ dont la transform\'ee de Fourier a un support dans $\{\xi\in\R^n:~\|\xi\|\le\s\}$. Alors il y a des fonctions $\Psi_j$, 
$1\le j\le n$, sur $\R^n\times\R^n$ qui appartiennent \`a l'espace $\fM_{\R^n\!,\R^n}$ de 
multiplicateurs de Schur et telles que
$$
f(x_1,\cdots,x_n)-f(y_1,\cdots,y_n)=\sum_{j=1}^n(x_j-y_j)\Psi_j(x_1,\cdots,x_n,y_1,\cdots,y_n),\quad x_j,~y_j\in\R,
$$
et
$$
\|\Psi_j\|_{\fM_{\R^n\!,\R^n}}\le\const\s\|f\|_{L^\be(\R^n)}.
$$
}

Ce r\'esultat implique que si $f$ appartient \`a l'espace de Besov $B_{\be,1}^1(\R^n)$ (voir \cite{Pee}), alors $f$ est \lb lipschitzienne op\'eratorielle, c'est-\`a-dire
$$
\|f(A_1,\cdots,A_n)-f(B_1,\cdots,B_n)\|\le\const\max_{1\le j\le n}\|A_j-B_j\|
$$
pour tous $n$-uplets $(A_1,\cdots,A_n)$ et $(B_1,\cdots,B_n)$ d'op\'erateurs auto-adjoints qui commutent.

Nous d\'emontrons aussi que si $f$ est une fonction h\"olderienne d'ordre $\a$ sur $\R^n$, alors $f$ est une fonction h\"olderienne op\'eratorielle d'ordre $\a$, c'est-\`a-dire
$$
\|f(A_1,\cdots,A_n)-f(B_1,\cdots,B_n)\|\le\const\max_{1\le j\le n}\|A_j-B_j\|^\a
$$
pour tous $n$-uplets $(A_1,\cdots,A_n)$ et $(B_1,\cdots,B_n)$ d'op\'erateurs autoadjoints qui commutent.

Nous obtenons aussi des analogues d'autres r\'esultats de \cite{Pe1}, \cite{Pe2}, \cite{AP2} et \cite{AP3} pour les fonctions d'$n$-uplets d'op\'erateurs auto-adjoint qui commutent (voir la version anglaise). 

\medskip

\begin{center}
------------------------------
\end{center}

\setcounter{section}{0}
\section{\bf Introduction}

\medskip

In this note we study the behavior of functions of perturbed tuples of commuting self-adjoint operators. We are going to find sharp estimates for $f(A_1,\cdots,A_n)-f(B_1,\cdots,B_n)$, where $(A_1,\cdots,A_n)$ and $(B_1,\cdots,B_n)$ are $n$-tuples of commuting self-adjoint operators and $f$ is a function on $\R^n$. Our results generalize the results of \cite{Pe1}, \cite{Pe2}, \cite{AP1}, \cite{AP2}, \cite{AP3}, \cite{AP4}, \cite{APPS1}, \cite{APPS2} for self-adjoint and normal operators.

Recall that a Lipschitz function $f$ on the real line $\R$ 
does not have satisfy the inequality
$$
\|f(A)-f(B)\|\le\const\|A-B\|
$$
for arbitrary self-adjoint operators $A$ and $B$ on Hilbert space, i.e., it does not have to be {\it operator Lipschitz}. This was proved in \cite{F}. Later it was shown in \cite{Pe1} and \cite{Pe2} that if $f$ is operator Lipschitz, then $f$ locally belongs to the Besov space $B_{1,1}^1(\R)$ (see \cite{Pee} for an introduction to Besov spaces) which also implies that Lipschitzness is not sufficient for operator Lipschitzness. On the other hand, it was proved in \cite{Pe1} and \cite{Pe2} that if $f$ belongs to the Besov space $B^1_{\infty,1}({\Bbb R})$, then $f$ is operator Lipschitz.

The situation changes dramatically if instead of the Lipschitz class, we consider the H\"older classes $\L_\a(\R)$, $0<\a<1$, of functions $f$ satisfying the inequality 
$|f(x)-f(y)|\le\const|x-y|^\a$, $x,\,y\in\R$. It was shown in \cite{AP1} and \cite{AP2} that a function $f$ in $\L_\a(\R)$ must be {\it operator H\"older of order $\a$}, i.e.,
$$
\|f(A)-f(B)\|\le\const\|A-B\|^\a.
$$
for arbitrary self-adjoint operators $A$ and $B$. Note that
the papers \cite{AP1} and \cite{AP2} also contain sharp estimates of $\|f(A)-f(B)\|$ for functions $f$ of class $\L_\o$ for arbitrary moduli of continuity $\o$.

It was also proved in \cite{AP1} and \cite{AP3} that if $f\in\L_\a$, $p>1$, and $A$ and $B$ are self-adjoint operators such that $A-B$ belongs to the Schatten--von Neumann class $\bS_p$, then $f(A)-f(B)\in\bS_{p/\a}$ and
$$
\|f(A)-f(B)\|_{\bS_{p/\a}}\le\const\|A-B\|_{\bS_p}^\a.
$$

Later in \cite{APPS1} and \cite{APPS2} the above results were generalized to the case of functions of normal operators. Note that the proofs given in \cite{Pe1}, \cite{Pe2}, \cite{AP1}, \cite{AP2}, and \cite{AP3} for self-adjoint operators do not work in the case of normal operators and a new approach was used in \cite{APPS1} and \cite{APPS2}.

In this paper we consider a more general problem of functions of $n$-tuples of commuting self-adjoint operators. The case $n=2$ corresponds to the case 
of normal operators. It turns out that the techniques used in \cite{APPS2} do not work for $n\ge3$. We offer in this note a new approach that works for all $n\ge1$.

We are going to use the technique of double operator integrals developed in \cite{BS1}, \cite{BS2}, and \cite{BS3}. Double operator integrals are expressions of the form
\bay
\label{doi}
\iint\limits_{\X_1\times\X_2}\Phi(s_1,s_2)\,dE_1(s_1)T\,dE_2(s_2),
\ey
where $E_1$ and $E_2$ are spectral measures on $\X_1$ and $\X_2$, $\Phi$ is a bounded measurable function on $\X_1\times\X_2$, and $T$ is an operator on Hilbert space. It was observed in \cite{BS1}, \cite{BS2}, and \cite{BS3} that the double operator integral \rf{doi} is well defined if $T\in\bS_2$ and determines an operator of class $\bS_2$. For certain $\Phi$, the transformer $T\mapsto \iint\Phi\,dE_1T\,dE_2$ maps the trace class $\bS_1$ into itself. If so, one can define by duality the integral \rf{doi} for all bounded operators $T$. Such functions $\Phi$ are called {\it Schur multipliers} (with respect to the spectral measures $E_1$ and $E_2$). We refer the reader to \cite{Pe1} for characterizations of Schur multipliers. 

If $\X_1$ and $\X_2$ are Borel subsets of Euclidean spaces, {\it we use the notation $\fM_{\X_1,\X_2}$ for the space of Borel functions $\Phi$ on $\X_1\times\X_2$ that are Schur multipliers for all Borel spectral measures $E_1$ and $E_2$ on $\X_1$ and $\X_2$}.

The proofs of the results of \cite{APPS2} for normal operators were based on the following formula: \begin{align*}
f(N_1)-f(N_2)=&\iint(\dg_yf)(z_1,z_2)\,dE_1(z_1)(B_1-B_2)\,dE_2(z_2)\\[.2cm]
&+
\iint(\dg_xf)(z_1,z_2)\,dE_1(z_1)(A_1-A_2)\,dE_2(z_2).
\end{align*}
Here $N_1$ and $N_2$ are normal operators with bounded difference $N_1-N_2$,  
$A_j=\re N_j$, $B_j=\im N_j$, $x_j=\re z_j$, $y_j=\im z_j$, $f$ is a bounded function on 
$\R^2$ whose Fourier transform has compact support, 
$$
(\dg_xf)(z_1,z_2)=\frac{f(x_1,y_2)-f(x_2,y_2)}{x_1-x_2},\quad\mbox{and}\quad
(\dg_yf)(z_1,z_2)=\frac{f(x_1,y_1)-f(x_1,y_2)}{y_1-y_2},\quad z_1,~z_2\in\C.
$$
It was shown in \cite{APPS2} that 
$\dg_xf$ and $\dg_yf$ belong to the space of Schur multipliers $\fM_{\R^2\!,\R^2}$.

However, 
in the case $n\ge3$ the situation is more complicated. Let $(A_1,A_2,A_3)$ and $(B_1,B_2,B_3)$ be triples of commuting self-adjoint operators. 
Suppose that $f$ is a bounded function on $\R^3$ whose Fourier transform has compact support. It can be shown that
\begin{align*}
f(A_1,A_2,A_3)&-f(B_1,B_2,B_3)=\iint(\dg_1f)(x,y)\,dE_1(x)(A_1-B_1)\,dE_2(y)\\[.2cm]
&+\iint (\dg_2f)(x,y)\,dE_1(x)(A_2-B_2)\,dE_2(y)
+\iint (\dg_3f)(x,y)\,dE_1(x)(A_3-B_3)\,dE_2(y),
\end{align*}
whenever the functions $\dg_1f$, $\dg_2f$, and $\dg_3f$ belong to the space of Schur multipliers $\fM_{\R^3\!,\R^3}$. Here 
$$
(\dg_1f)(x,y)=\frac{f(x_1,x_2,x_3)-f(y_1,x_2,x_3)}{x_1-y_1},\quad
(\dg_2f)(x,y)=\frac{f(y_1,x_2,x_3)-f(y_1,y_2,x_3)}{x_2-y_2},
$$
$$
(\dg_3f)(x,y)=\frac{f(y_1,y_2,x_3)-f(y_1,y_2,y_3)}{x_3-y_3},
\quad x=(x_1,x_2,x_3),\quad y=(y_1,y_2,y_3).
$$
The methods of \cite{APPS2} allow us to prove that $\dg_1f$ and $\dg_3f$ do 
belong to the space of Schur multipliers $\fM_{\R^3\!,\R^3}$. However, as the next result shows, the function $\dg_2f$ does not have to be in $\fM_{\R^3\!,\R^3}$.

\begin{thm}
\label{kontr}
Suppose that $g$ is a bounded function on $\R$ such that the Fourier transform of $g$ has compact support and is not a measure. Let $f$ be the function on $\R^3$ defined by
$$
f(x_1,x_2,x_3)=g(x_1-x_3)\sin x_2.
$$ 
Then $\dg_2f\not\in\fM_{\R^3\!,\R^3}$.
\end{thm}

Note that it is easy to construct such a function $g$, e.g., $g(x)=\int_0^xt^{-1}\sin t\,dt$.

In \S\,2 we show that in the case $n\ge3$ it is possible to represent 
$f(A_1,\cdots,A_n)-f(B_1,\cdots,B_n)$ in terms of double operator integrals in a different way. Using such a representation, we obtain in \S\,3 and \S\,4 analogs of the above results  in the case of $n$-tuples of commuting self-adjoint operators.

\medskip

\section{\bf An integral representation}

\medskip

The integral representation for $f(A_1,\cdots,A_n)-f(B_1,\cdots,B_n)$ is based on the following result:

\begin{thm}
\label{fsi}
Let $\s>0$ and let $f$ be a function in $L^\be(\R^n)$ whose Fourier transform is supported on 
$\{\xi\in\R^n:~\|\xi\|\le\s\}$. Then there exist functions $\Psi_j$ in $\fM_{\R^n\!,\R^n}$, $1\le j\le n$, such that
\bay
\label{psij}
f(x_1,\cdots,x_n)-f(y_1,\cdots,y_n)=\sum_{j=1}^n(x_j-y_j)\Psi_j(x_1,\cdots,x_n,y_1,\cdots,y_n),\quad x_j,~y_j\in\R,
\ey
and $\|\Psi_j\|_{\fM_{\R^n\!,\R^n}}\le\const\s\|f\|_{L^\be(\R^n)}$.
\end{thm}

\newcommand\mC{\mathcal{C}}
\newcommand\mQ{\mathcal{Q}}
\newcommand\mR{\mathcal{R}}

We are going to derive Schur multiplier estimates from the following lemma.

\begin{lem}
\label{Fou}
Let $\mC=\mQ\times\mR$ be a cube in $\R^{2n}$ of sidelength $L$ and let $\Psi$ be a $C^\be$ function on $\frac{3}{2}\mC$. Then $\Psi\big|\mC\in\fM_{\mQ,\mR}$ and
$$
\|\Psi\|_{\fM_{\mQ,\mR}}\le\const\max\Big\{L^{|\a|}\max_{a\in\frac32\mC}\big|(D^\a\Psi)(a)\big|:~|\a|\le2n+2\Big\}.
$$
\end{lem}

The lemma can be proved by expanding $\Psi$ in the Fourier series.

{\bf Sketch of the proof of Theorem \ref{fsi}.} By rescaling, we may assume that 
$\|f\|_{L^\be}\le1$ and $\s=1$.

We consider the lattice of dyadic cubes in $\R^{2n}=\R^n\times\R^n$, i.e., the cubes whose sides are intervals of the form 
$\big[j2^k,(j+1)2^k\big)$, $j,\,k\in\Z$. We say that a dyadic cube $\mC$ in $\R^{2n}=\R^n\times\R^n$ is {\it admissible} if either its sidelength ${\rm L}(\mC)$ is equal to $1$ or ${\rm L}(\mC)>1$ and the interior of the cube $2\mC$, i.e., the cube centered at the center of $\mC$ with sidelength $2{\rm L}(\mC)$, does not intersect the diagonal $\{(x,x):~x\in\R^n\}$. An admissible cube is called {\it maximal} if it is not a proper subset of another admissible cube. It is easy to see that the maximal admissible cubes are disjoint and cover $\R^{2n}$. It can also easily be verified that if $\mQ$ is a dyadic cube in $\R^n$, then there can be at most $6^n$ dyadic cubes $\mR$ in $\R^n$ such that $\mQ\times\mR$ is a maximal admissible cube.
For $l=2^m$, we denote by $\cd_l$ the set of maximal dyadic cube of sidelength $l$.

It follows that if $\O$ is a function on $\R^n\times\R^n$ that is supported on
$\bigcup_{\mC\in\cd_l}\mC$, then 
$$
\|\O\|_{\fM_{\R^n\!,\R^n}}\le6^n\sup_{\mC\in\cd_l}\|\chi_{_\mC}\O\|_{\fM_{\R^n\!,\R^n}}.
$$

We have to define $\Psi_j$ on each maximal admissible cube. Suppose that $\mC\in\cd_1$. We put
$$
\Psi_j(x,y)=\int_0^1(D_jf)((1-t)x+ty)\,dt,\quad(x,y)\in\mC=\mQ\times\mR,
$$
where $D_jf$ is the $j$th partial derivative of $f$. It follows from Lemma \ref{Fou} that 
$\|\Psi_j\|_{\fM_{\mQ,\mR}}\le\const$.

Suppose now that $l=2^m>1$ and $\mC=\mQ\times\mR\in\cd_l$. 
Let $\o$ be a $C^\be$ nonnegative even function on $\R$ such that 
$\o(t)=0$ for $t\in[-\frac12,\frac12]$, and $\o(t)=1$ for $t\not\in[-1,1]$. We put 
$\Phi_j(x,y)=\o((x_j-y_j)/l)$,
$\Phi=\sum_{j=1}^n\Phi_j$, and define the functions $\Xi_j$, $1\le j\le n$, by
$$
\Xi_j(x,y)=\left\{\begin{array}{ll}
\frac{1}{x_j-y_j}\cdot\frac{\Phi_j(x,y)}{\Phi(x,y)},&x_j\ne y_j,\\[.2cm]
0,&x_j=y_j.\end{array}\right.
$$
It follows easily from Lemma \ref{Fou} that $\|\Xi_j\|_{\fM_{\mQ,\mR}}\le\const 2^{-m}$. We put now 
$$
\Psi_j(x,y)=(f(x)-f(y))\Xi_j(x,y),\quad(x,y)\in\mC.
$$
Clearly, \rf{psij} holds for $(x,y)\in\mC$ and 
$\|\Psi_j\|_{\fM_{\mQ,\mR}}\le\const 2^{-m}$. The functions $\Psi_j$ are now defined on $\R^n\times\R^n$ and 
$
\|\Psi_j\|_{\fM_{\R^n,\,\R^n}}\le\const\sum_{m\ge0}2^{-m}.
$
This implies the result. $\bl$


\begin{thm}
\label{intpr}
Let $f$ be a function satisfying the hypotheses of Theorem {\em\ref{fsi}} and let $\Psi_j$, $1\le j\le n$, be functions in $\fM_{\R^n\!,\R^n}$ satisfying {\em\rf{psij}}. Suppose that
$(A_1,\cdots,A_n)$ and $(B_1,\cdots,B_n)$ are $n$-tuples of commuting self-adjoint operators such that the operators $A_j-B_j$ are bounded, $1\le j\le n$.
Then
$$
f(A_1,\cdots,A_n)-f(B_1,\cdots,B_n)=\sum_{j=1}^n\,\,\,\,~\iint\limits_{\R^n\times\R^n}\Psi_j(x,y)\,dE_A(x)(A_j-B_j)\,dE_B(y)
$$
and
$
\|f(A_1,\cdots,A_n)-f(B_1,\cdots,B_n)\|\le\const\s\|f\|_{L^\be(\R^n)}
\max_{1\le j\le n}\|A_j-B_j\|.
$
\end{thm}

\medskip

\section{\bf Operator norm estimates}

\medskip

In this section we obtain operator norm estimates for $f(A_1,\cdots,A_n)-f(B_1,\cdots,B_n)$, where \lb$(A_1,\cdots,A_n)$ and $(B_1,\cdots,B_n)$ are $n$-tuples of commuting self-adjoint operators.

A function $f$ on $\R^n$ is called {\it operator Lipschitz} if 
$$
\|f(A_1,\cdots,A_n)-f(B_1,\cdots,B_n)\|\le\const\max_{1\le j\le n}\|A_j-B_j\|
$$
for all $n$-tuples of commuting self-adjoint operators $(A_1,\cdots,A_n)$ and $(B_1,\cdots,B_n)$.

The following theorem can be deduced easily from Theorem \ref{intpr}. 

\begin{thm}
\label{besov}
Let $f$ be a function in the Besov space $B_{\be,1}^1(\R^n)$. Then $f$ is operator Lipschitz.
\end{thm}

For $\a\in(0,1)$, we define the H\"older class $\L_\a(\R^n)$ of functions $f$ on $\R^n$ such that
$$
|f(x)-f(y)|\le\const\|x-y\|_{\R^n}^\a,\quad x,~y\in\R^n.
$$
For a modulus of continuity $\o$, the space $\L_\o(\R^n)$ consists of functions $f$ on $\R^n$ such that 
$$
|f(x)-f(y)|\le\const\o\big(\|x-y\|_{\R^n}\big),\quad x,~y\in\R^n.
$$

The following results are analogs of the corresponding results of \cite{AP1} and \cite{AP2}
in the case $n=1$. The proofs of Theorems \ref{alpha} and \ref{omega} are based on Theorem \ref{intpr} and use the same methods as in \cite{AP2}.

\begin{thm}
\label{alpha}
Let $\a\in(0,1)$ and let $f\in\L_\a(\R^n)$. Then $f$ is operator H\"older of order $\a$, i.e.,
$$
\|f(A_1,\cdots,A_n)-f(B_1,\cdots,B_n)\|\le\const\max_{1\le j\le n}\|A_j-B_j\|^\a
$$
for all $n$-tuples of commuting self-adjoint operators $(A_1,\cdots,A_n)$ and $(B_1,\cdots,B_n)$.
\end{thm}

\begin{thm}
\label{omega}
Let $\o$ be a modulus of continuity and let $f\in\L_\o(\R^n)$. Then 
$$
\|f(A_1,\cdots,A_n)-f(B_1,\cdots,B_n)\|
\le\const\o_*\left(\max_{1\le j\le n}\|A_j-B_j\|\right)
$$
for all $n$-tuples of commuting self-adjoint operators $(A_1,\cdots,A_n)$ and $(B_1,\cdots,B_n)$, where
$$
\o_*(\d)\df\d\int_\d^\be\frac{\o(t)}{t^2}\,dt,\quad\d>0.
$$
\end{thm}

\medskip

\section{\bf Schatten--von Neumann norm estimates}

\medskip

In this section we obtain estimates in $\bS_p$ norms. 

\begin{thm}
\label{S1}
Let $f$ be a function in the Besov space $B_{\be,1}^1(\R^n)$. Suppose that 
$(A_1,\cdots,A_n)$ and $(B_1,\cdots,B_n)$ are $n$-tuples of commuting self-adjoint operators such that $A_j-B_j\in\bS_1$. Then $f(A_1,\cdots,A_n)-f(B_1,\cdots,B_n)\in\bS_1$ and 
$$
\|f(A_1,\cdots,A_n)-f(B_1,\cdots,B_n)\|_{\bS_1}\le\const\|f\|_{B_{\be,1}^1(\R^n)}
\max_{1\le j\le n}\|A_j-B_j\|_{\bS_1}.
$$
\end{thm}

\begin{thm}
\label{Spa}
Let $f\in\L_\a(\R^n)$ and let $p>1$. Suppose that 
$(A_1,\cdots,A_n)$ and $(B_1,\cdots,B_n)$ are $n$-tuples of commuting self-adjoint operators such that $A_j-B_j\in\bS_p$. Then 
$f(A_1,\cdots,A_n)-f(B_1,\cdots,B_n)\in\bS_{p/\a}$ and 
$$
\|f(A_1,\cdots,A_n)-f(B_1,\cdots,B_n)\|_{\bS_{p/\a}}\le\const\|f\|_{\L_\a(\R^n)}
\max_{1\le j\le n}\|A_j-B_j\|^\a_{\bS_p}.
$$
\end{thm}

Note that the conclusion of Theorem \ref{Spa} does not hold in the case $p=1$ even if $n=1$, see \cite{AP3}.

\begin{thm}
\label{Ba}
Let $f$ be a function in the Besov space $B_{\be,1}^\a(\R^n)$. Suppose that 
$(A_1,\cdots,A_n)$ and $(B_1,\cdots,B_n)$ are $n$-tuples of commuting self-adjoint operators such that $A_j-B_j\in\bS_1$. Then 
$f(A_1,\cdots,A_n)-f(B_1,\cdots,B_n)\in\bS_{1/\a}$ and 
$$
\|f(A_1,\cdots,A_n)-f(B_1,\cdots,B_n)\|_{\bS_{1/\a}}\le\const\|f\|_{B_{\be,1}^\a(\R^n)}
\max_{1\le j\le n}\|A_j-B_j\|^\a_{\bS_1}.
$$
\end{thm}

The proofs of the above theorems are based on Theorem \ref{intpr} and use the methods of \cite{AP3}.

Note that in \cite{AP3} more general results for other operator ideals were obtained in the case $n=1$. Those results can also be generalized to the case of arbitrary $n\ge1$.

We would like to mention the paper \cite{KPSS} on Lipschitz estimates in the norm of $\bS_p$,
$1<p<\be$, for functions of commuting tuples of self-adjoit operators.

\end{document}

%% file: CRinput.tex
\newcommand{\lb}{\linebreak}

\renewcommand{\a}{\alpha}

\renewcommand{\d}{\delta}

\newcommand{\s}{\sigma}

\renewcommand{\o}{\omega}

\renewcommand{\L}{\Lambda}

\renewcommand{\O}{\Omega}

\newcommand{\cd}{{\mathcal D}}

\newcommand{\X}{{\mathcal X}}

\newcommand{\C}{{\Bbb C}}

\newcommand{\R}{{\Bbb R}}
\newcommand{\Z}{{\Bbb Z}}

\newcommand{\bS}{{\boldsymbol S}}

\newcommand{\rf}[1]{(\ref{#1})}

\newcommand{\df}{\stackrel{\mathrm{def}}{=}}

\newcommand{\re}{\operatorname{Re}}

\newcommand{\const}{\operatorname{const}}

\newcommand{\eeq}{\end{equation}}
\newcommand{\beq}{\begin{equation}}
\newcommand{\bay}{\begin{eqnarray}}
\newcommand{\ba}{\begin{align*}}
\newcommand{\ea}{\end{align*}}
\newcommand{\ey}{\end{eqnarray}}
\newcommand{\bey}{\begin{eqnarray*}}
\newcommand{\eey}{\end{eqnarray*}}

\newcommand{\be}{\infty}

\newcommand{\bl}{\blacksquare}

\newcommand{\im}{\operatorname{Im}}
\renewcommand{\re}{\operatorname{Re}}

\newtheorem{thm}{\hspace{\parindent}Theorem}[section]

\newtheorem{lem}[thm]{\hspace{\parindent}Lemma}